\documentclass[12pt,leqno]{amsart}
\usepackage{amscd}
\usepackage{color}
\usepackage{amssymb}
\usepackage[matrix,arrow]{xy}

\newtheorem{theorem}{Theorem}[section]
\newtheorem{proposition}[theorem]{Proposition}
\newtheorem{lemma}[theorem]{Lemma}

\newtheorem{corollary}[theorem]{Corollary}

\newtheorem*{mainthm}{Theorem}

\theoremstyle{remark}
\newtheorem{remark}[theorem]{Remark}

\newtheorem{definition}[theorem]{Definition}

\newcommand{\U}{\mathcal{U}}
\newcommand{\V}{\mathcal{V}}
\newcommand{\R}{\mathbb{R}}
\newcommand{\UE}{uniformly embeddable}
\newcommand{\UUE}{equi-uniformly embeddable }
\newcommand{\eps}{\varepsilon}

\newcounter{ictr}

\newenvironment{ilist}{\begin{list}
                         {\textup{(\roman{ictr})}}
                         {\usecounter{ictr}
                          \setlength{\leftmargin}{0.6truein}
                          \setlength{\itemsep}{0.0truein}
                          \setlength{\labelwidth}{0.3truein}}}
                      {\end{list}}

\newcounter{actr}
\newenvironment{alist}{\begin{list}
                         {\textup{(\alph{actr})}}
                         {\usecounter{actr}
                          \setlength{\leftmargin}{0.6truein}
                          \setlength{\itemsep}{0.0truein}
                          \setlength{\labelwidth}{0.3truein}}}
                      {\end{list}}

\newcommand{\Hawaii}{Hawai\kern.05em`\kern.05em\relax i}
\newcommand{\Manoa}{M\=anoa}

\begin{document}
\title{ Uniform embeddability of relatively hyperbolic groups}

\author{Marius Dadarlat}
  \address{ Department of Mathematics,
Purdue University, 150 N. University Street, West Lafayette, IN
47907-2067}
 \email{mdd@math.purdue.edu}
\author{Erik Guentner}
\address{Mathematics Department, University of \Hawaii, \Manoa,
  2565 McCarthy Mall, Honolulu, HI 96822}
\email{erik@math.hawaii.edu}

\thanks{The first author was supported in part by NSF Grant
\#DMS-0200601.  The second author was supported in part by NSF
Grant \#DMS-0349367}
\date{\today}

\begin{abstract}
Let $\Gamma$ be a finitely generated group which is hyperbolic
relative to a finite family $\{H_1,\dots,H_n\}$ of subgroups. We
prove that $\Gamma$ is uniformly embeddable in a Hilbert space if
and only if each subgroup $H_i$ is uniformly embeddable in a
Hilbert space.
\end{abstract}

\maketitle


\section{Introduction}

Gromov introduced the notion of uniform embeddability (in Hilbert
space), and suggested it should be relevant to the Novikov
Conjecture \cite{Grom,NCITR}.  Subsequently Yu proved the Coarse
Baum-Connes Conjecture for bounded geometry discrete metric spaces
that are uniformly embeddable; applying a descent principle the
Novikov Conjecture followed for groups that, when equipped with a
word length metric, are uniformly embeddable \cite{Yu:BC}.  (A
condition on finiteness of the classifying space was later removed
\cite{Skan-Tu-Yu}.)

The notion of $C^*$-exactness of discrete groups was introduced by
Kirchberg.  It has been extensively studied as a functional
analytic property of groups, and developed by many authors.  In
particular, Ozawa gave a characterization of $C^*$-exact groups,
from which it directly follows that a $C^*$-exact group is
uniformly embeddable \cite{Ozawa-char} (see \cite{Gue-Kam} for the
link to uniform embeddability).  There is at present no known
example of a group that is uniformly embeddable but not
$C^*$-exact. Analogous statements for metric spaces involve
Property $A$ of Yu \cite{Yu:BC}, which is equivalent to
$C^*$-exactness for discrete groups.

The classes of metric spaces (and groups) that are uniformly
embeddable, or have Property $A$ (are $C^*$-exact) are the subject
of intense study.  In this note we introduce a `gluing' technique
for proving uniform embeddability: starting from the assumption
that a space is covered in an appropriate way by uniformly
embeddable sets we conclude that the space itself is uniformly
embeddable.  Thus, the individual uniform embeddings of the pieces
are `glued' to give a uniform embedding of the whole.  A parallel
technique is introduced for spaces with Property $A$; it applies
to $C^*$-exact groups.  The most primitive gluing result is
summarized in the following theorem (compare
Theorem~\ref{UE-gluing}):

\begin{mainthm}
  Let $X$ be a metric space.  Assume that for all $\lambda>0$ there
  exists a partition of unity $(\varphi_i)_{i\in I}$ on $X$ such that
  \begin{ilist}
    \item the associated $\Phi:X\to l^1(I)$
       \textup{(}defined by $\Phi(x)(i)=\varphi_i(x)$\textup{)}
        is Lipschitz, and
    \item the subspaces $(\mathrm{supp}(\varphi_i))_{i\in I}$ are
          `equi' uniformly embeddable.
  \end{ilist}
  Then $X$ is uniformly embeddable.
\end{mainthm}

Our gluing technique is inspired by the work of Bell and
Dranishnikov on spaces of finite asymptotic dimension
\cite{Bell-Dran,Roe}, and by the work of Bell on Property $A$
\cite{Bell}.  Our results differ from these works in several ways:
first, we treat uniformly embeddable spaces, and allow spaces of
unbounded geometry; second, we significantly relax the condition
on the `parameter space'.  With these refinements, the gluing
technique is extremely versatile --- it allows us to give a
conceptual treatment of all known permanence properties of the
classes of uniformly embeddable spaces and $C^*$-exact groups
(see, for example, \cite{Dad-Guen,Dyk,Kirch-Wass}).

In this note we concentrate on the basics of our gluing technique;
we plan to further develop the technique and present the above
mentioned applications elsewhere.  We also describe one
application of gluing, proving the following new permanence
property of the class of uniformly embeddable groups:

\begin{mainthm}
Let $\Gamma$ be a finitely generated group which is hyperbolic
relative to a finite family $\{H_1,\dots,H_n\}$ of subgroups. Then
$\Gamma$ is uniformly embeddable if and only if each subgroup
$H_i$ is uniformly embeddable.
\end{mainthm}

A parallel statement concerning $C^*$-exactness was recently
obtained by Ozawa \cite{Ozawa}.  Our methods allow us to recover
this result; indeed, the $C^*$-exactness result may be extracted
from the combined work of Bell \cite{Bell} and Osin \cite{Osin}.
Our approach to relative hyperbolicity, which is based on the work
of Osin, is quite different from that taken by Ozawa however, who,
given amenable $H_i$-spaces explicitly constructs an amenable
$\Gamma$-space.

These results should be compared to a recent result of Osin
\cite{Osin}: in the situation described in the theorem, $\Gamma$
has finite asymptotic dimension if and only the $H_i$ have finite
asymptotic dimension.

\section{Preliminaries}

Let $X$ and $Y$ be metric spaces, with metrics $d_X$ and $d_Y$,
respectively.  A function $F: X \to Y$ is a {\it uniform
embedding\/} if there exist non-decreasing functions $\rho_\pm
:\R_+\to \R_+$ such that $\lim_{\,t \to \infty}\rho_\pm(t)=\infty$
and such that
\begin{equation}
\label{eq:ue}
  \rho_-(d_X(x,x'))\leq  d_Y(F(x),F(x')) \leq \rho_+(d_X(x,x')),
      \quad\text{for all $x$, $x'\in X$.}
\end{equation}
The space $X$ is (Hilbert space) {\it uniformly embeddable\/}  if
there exists a uniform embedding $F$ of $X$ into a (real) Hilbert space
$\mathcal{H}$.
\begin{proposition}\label{char-UE}
  Let $X$ be a metric space.  Then $X$ is uniformly embeddable if and
  only if for every $R >0$ and $\varepsilon>0$ there exists a Hilbert space
  valued map $\xi:X \to \mathcal{H}$, $(\xi_x)_{x \in X}$, such that
  $\|\xi_x\|=1$, for all $x \in X$ and such that
\begin{ilist}
  \item $d(x,x') \leq R \Rightarrow \|\xi_x-\xi_{x'}\| \leq \eps$
  \item $\lim_{S \to \infty} \sup \{ |\langle\xi_x,\xi_{x'}\rangle|:
      d(x,x') \geq S,\, x,x'\in X\} =0 .$\qed
\end{ilist}
\end{proposition}
\begin{definition} \label{UUE}
A family of metric spaces $(X_i,d_i)$ is called \UUE\ if there is
a family of Hilbert space valued maps   $F_i: X_i \to H_i$ and
non-decreasing functions $\rho_\pm :\R_+\to \R_+$ with $\lim_{\,t
\to \infty}\rho_\pm(t)=\infty$ such that
\begin{equation}
\label{eq:uue}
  \rho_-(d_i(x,y))\leq  \|F_i(x)-F_i(y)\| \leq \rho_+(d_i(x,y)),
      \quad\text{for all $i$ and all $x$, $y\in X_i$.}
\end{equation}
\end{definition}
The proof of Proposition~\ref{char-UE} given in \cite{Dad-Guen}
for arbitrary metric spaces also shows the following
characterization of equi-uniform embeddability
 for families of metric spaces.
\begin{proposition}\label{char-UUE}
     A family $(X_i)_{i\in I}$ of metric spaces  is equi-uniformly embeddable if and
  only if for every $R >0$ and $\varepsilon>0$ there exists a family
  $(\xi_i)_{i \in I}$ of Hilbert space
  valued maps $\xi_i:X_i \to \mathcal{H}$,   such that
  $\|\xi_i(x)\|=1$, for all $x \in X_i$,  and such that
\begin{ilist}
  \item $\forall \,i \in I$ $\forall x,x' \in X_i$
  $d(x,x') \leq R \Rightarrow \|\xi_i(x)-\xi_i(x')\| \leq \eps$
  \item $\lim_{S \to \infty} \sup_{i \in I}\sup \{ |\langle\xi_i(x),\xi_i(x')\rangle|:
      d(x,x') \geq S,\, x,x'\in X_i\} =0 .$ \qed
\end{ilist}
\end{proposition}

\begin{remark}\label{equi-UE}
Let $X$ be a metric space which is uniformly embeddable.
Then any family $(X_i)_{i\in I}$ of subspaces of
$X$ is equi-uniformly embeddable.
\end{remark}

Property $A$ is a condition on metric spaces introduced by Yu
\cite{Yu:BC}.  We do not recall the definition of Property $A$
here; rather, we work with the following characterization of
Property $A$ obtained by Tu \cite{Tu:notes}.

\begin{proposition}[\cite{Tu:notes}]\label{EXACT}
A discrete metric space $X$ with bounded geometry has Property A
if and only if for every $R>0$ and $\varepsilon>0$ there exist a
function $\xi:X\to \ell^1(X)$  and a number $S>0 $ such that for
all $x$, $x'\in X$ we have $\| \xi_x\|=1$, and
\begin{ilist}
  \item $d(x,x') \leq R \Rightarrow \|\xi_x-\xi_{x'}\| \leq \eps$
  \item $\mathrm{supp}\, \xi_x\subset B(x,S)$.
\end{ilist}
\textup{(}Moreover one can arrange that $\xi$ is
nonnegative.\textup{)}

Equivalently, for every $R>0$ and $\eps>0$ there exists an $S>0$
and a Hilbert space valued function $\xi:X\to \mathcal{H}$ such
that for all $x$, $x'\in X$ we have $\|\xi_x\|=1$,
\textup{(}i\textup{)} as above
and
\begin{ilist}{\setcounter{ictr}{2}}
 \item $\exists$
$S>0$ such that $d(x,x') \geq S \Rightarrow
          \langle\, \xi_x, \xi_{x'} \,\rangle= 0$.\qed
\end{ilist}
\end{proposition}

\begin{remark}
The existence of $\xi$ satisfying the conditions (i), (ii)
(respectively (iii)) of Proposition~\ref{EXACT} is a consequence
of Property $A$ for {\it arbitrary\/} metric spaces as shown in
\cite{Tu:notes}. The bounded geometry condition is needed only for
the reverse implication.
\end{remark}

Let $X$ be a set. A {\it partition of unity\/} on $X$ is a family
of maps $(\varphi_i)_{i\in I}$, with $\varphi_i:X \to [0,1]$, and
such that $\sum_{i\in I}\varphi_i (x)=1$ for all $x \in X$. If $x
\in X$ we do not require that the set $\{i \in I:\varphi_i(x)\neq
0\}$ be finite, although that will be the case in most of our
examples. We say that $(\varphi_i)_{i\in I}$ is {\it subordinated
to a cover\/} $\U=(U_i)_{i \in I}$ of $X$ if each $\varphi_i$
vanishes outside $U_i$. Sometimes a partition of unity
subordinated to a cover $\U$ will be denoted by $(\varphi_U)_{U
\in \U}$.

\begin{definition}\label{def-exactness}
A metric space $X$ is \emph{exact} if for all $R>0$ and
$\varepsilon>0$ there is a partition of unity $(\varphi_i)_{i\in
I}$ on $X$ subordinated to a cover $\U=(U_i)_{i \in I}$ and such
that
\begin{ilist}
  \item $\forall x,y\in X$ with $d(x,y)\leq R$,
    $\sum _{i\in I}|\varphi_i (x)-\varphi_i (y)|\leq \varepsilon.$
  \item the cover $\U= (U_i)_{i\in I}$ is a uniformly bounded, i.e.
    $\sup_{i \in I} \mathrm{diam}(U_i) < \infty$.
\end{ilist}
\end{definition}

\begin{definition}\label{def-unif-exactness}
A family of metric spaces  $(X_i)_{i \in I}$ is \emph{equi-exact}
if for all $R>0$ and $\varepsilon>0$ and for every $i \in I$ there
is a partition of unity $(\psi_i^j)_{j\in\, J_i}$ on $X_i$
subordinated to a cover $\U_i=(U_i^j)_{j\in \,J_i}$ of $X_i$ and
such that
\begin{ilist}
  \item $\forall i \in I$, $\forall x,y\in X_i$
       with $d(x,y)\leq R$,
       $\sum _{j\in \,J_i}|\psi_i^j (x)-\psi_i^j(y)|
             \leq  \varepsilon$.
  \item the family $(U_i^j)_{i\in I,\, j\in \,J_i}$ is uniformly
        bounded, i.e.
       $\sup_{i\in I,\, j\in \,J_i} \mathrm{diam}(U_i^j) < \infty$.
\end{ilist}
\end{definition}

\begin{remark}\label{equi-exact}
Let $X$ be an exact metric space. Then any family $(X_i)_{i\in I}$
of subspaces of $X$ is equi-exact.
\end{remark}

\begin{proposition} Let $X$ be a metric space.
\begin{alist}
  \item If $X$ has  Property A  then $X$ is exact.
  \item If $X$ is discrete and has bounded geometry then $X$ is exact
        if and only if it has Property A.
  \item If $X$ is exact then $X$ is uniformly embeddable.
\end{alist}
\end{proposition}



\begin{proof}
For (a), assume that $X$ has Property A. Let $R>0$ and $\eps>0$ be given.
Obtain
$\xi:X \to \ell^1(X)$ (which we assume to be non-negative)
and $S>0$ as in Proposition~\ref{EXACT}. Define
$(\varphi_z)_{z \in X}$ by $\varphi_z(x)=\xi_x(z)$. If $U_z=\{x\in
X:\varphi_z(x)>0\}$ then $U_z\subset B(z,S)$ since $\mathrm{supp}\, \xi_x
\subset B(x,S)$. It is clear that $\sum_{z \in
X}\,\varphi_z(x)=\|\xi_x\|=1$ and $$\sum _{z\in X}|\varphi_z
(x)-\varphi_z (y)|=\|\xi_x-\xi_y\|\leq
  \varepsilon$$
  if $d(x,y)\leq R$.

For (b), assume that $X$ is exact. It suffices to find a Hilbert
space valued function on $X$ satisfying the conditions from the
second part of Proposition~\ref{EXACT}. Let $R>0$ and $\eps>0$ be
given.  Let $(\varphi_i)_{i\in I}$ be as in
Definition~\ref{def-exactness}. Define $\xi:X \to \ell^2(I)$ by
$\xi_x(i)=\varphi_i(x)^{1/2}$. Then $\|\xi_x \|^2= \sum_{i \in
I}\,\varphi_i(x) =1$. Using the inequality $|a^{1/2} - b^{1/2}|^2
\leq |a-b|$ we see that
$$\|\xi_x-\xi_y\|^2=\sum_{i\in I }|\varphi_i(x)^{1/2}-\varphi_i(y)^{1/2}|^2\leq
    \sum _{i\in I}|\varphi_i (x)-\varphi_i (y)|\leq
  \varepsilon$$ if $d(x,y)\leq R$.
  Finally we note that
  $$\langle \xi_x,\xi_y\rangle=\sum _{i\in I}\,\varphi_i (x)^{1/2} \varphi_i
  (y)^{1/2}=0$$ whenever $d(x,y) >\sup_{i \in I} \mathrm{diam}(U_i)$.

The proof of (c) is similar to that of (b), but one applies
  Proposition~\ref{char-UE}.
\end{proof}


\begin{remark}\label{geometric}
If $f:X \to Y$ is a uniform embedding of metric spaces and $Y$ is
exact then $X$ is exact. Therefore exactness of a metric space is a
coarse invariant. This remark is generalized in
Corollary~\ref{fibrewise-general}.
\end{remark}

\section{Gluing spaces using partitions of unity}

Let $X$ be a metric space and let $\U=(U_i)_{i\in I}$ be a cover
of $X$. Denote by $\U_R$ the cover obtained by enlarging the sets
in $\U$ by taking their $R$-closed neighborhoods:
\begin{equation*}
  \U_R = \{\, U_i(R) \colon i \in I \,\}, \quad
     U_i(R) = \{\, x\in X \colon d(x,U_i)\leq R \,\}.
\end{equation*}
One verifies immediately that if the family $\U$ (with the metric
structure induced from $X$) is equi-uniformly embeddable
 (or equi-exact) then so is $\U_R$.

\begin{theorem} \label{exact-gluing}
Let $X$ be a metric space.  Suppose that for all $R>0$ and
$\varepsilon>0$ there is a partition of unity
$(\varphi_i)_{i\in I}$ on $X$ such that
\begin{ilist}
  \item $\forall x,y\in X$ with $d(x,y)\leq R$,
     $\sum _{i\in I}|\varphi_i (x)-\varphi_i (y)|\leq \varepsilon$,
     and
   \item    $(\varphi_i)_{i\in I}$ is subordinated to an equi-exact
    cover $(U_i)_{i\in I}$ of $X$.
\end{ilist}
Then $X$ is exact.
\end{theorem}

\begin{proof}
Let $R>0$ and $\epsilon >0$ be given. We are going to construct a
Hilbert space valued function $\eta:A \to \mathcal{H}$ satisfying
the hypotheses of Proposition~\ref{char-UE}.
 By assumption there is a cover
$\U=(U_i)_{i\in I}$ of $X$ which is equi-exact and there is a
partition of unity $(\varphi_i)_{i\in I}$ subordinated to $\U$
such that
  $\forall x,y\in X$ with $d(x,y)\leq R$,
\begin{equation*}
\sum_{i\in I}|\varphi_i (x)-\varphi_i (y)|\leq
  \varepsilon/2.
\end{equation*}
Since $\U$ is equi-exact so is $\U_R=\{U_i(R):i \in I\}$.
Therefore for each $U_i \in \U$ there is a cover $\V_i
=(V_{i}^j)_{j \in\, J_i}$ of $U_i(R)$ such that the cover $\{\,
V_i^j \colon i\in I,\, j\in \,J_i \,\} $ of $X$ is uniformly
bounded. Moreover for each $U_i \in \U$  there is a partition of
unity $(\psi_{i}^j)$ on $U_i(R)$ subordinated to $\V_i$ such that
$\forall x,y\in U_i(R)$ with $d(x,y)\leq R$,
\begin{equation*}
\sum _{j\in\, J_i}|\psi_{i}^j(x)-\psi_{i}^j(y)|\leq
  \varepsilon/2.
\end{equation*}
It is useful to extend $\psi_{i}^j$ to $X$ by setting it equal to
zero outside $U_i(R)$.  Define $\theta_{i,j}=\varphi_i\psi_{i}^j$.
Then $(\theta_{i,j})$ is a partition of unity on $X$ subordinated
to a uniformly bounded cover. Moreover,
\begin{equation*}
\sum_{i,j}|\theta_{i,j}(x)-\theta_{i,j}(y)|\leq
  \sum_{i}\varphi_i(x)\sum_{j}|\psi_{i}^j(x)-\psi_{i}^j(y)|+
  \sum_{i} |\varphi_{i} (x)-\varphi_{i}(y)|\sum_{j}  \psi_{i}^j(y).
\end{equation*}
Assume now that $d(x,y)\leq R$. If $\varphi_{i} (x)\neq 0$ then $x
\in U_i$ hence $y \in U_i(R)$ as $d(x,y)\leq R$. Therefore
\begin{equation*}
\sum_{i}\varphi_i(x)\sum_{j}|\psi_{i}^j(x)-\psi_{i}^j(y)|\leq
\varepsilon/2.
\end{equation*}
Since $\sum_j\psi_{i}^j(y) $ equals $1$ for $y \in U_i(R)$ and $0$
for $y \notin U_i(R)$,
\begin{equation*}
\sum_{i} |\varphi_{i} (x)-\varphi_{i}(y)|\sum_{j}  \psi_{i}^j(y)
  \leq \sum_{i} |\varphi_{i} (x)-\varphi_{i}(y)|\leq
  \varepsilon/2.
\end{equation*}
  Combining the above estimates we obtain that
  $$\sum_{i,j}|\theta_{i,j}(x)-\theta_{i,j}(y)|\leq \varepsilon$$
  whenever $x,y \in X$ and $d(x,y)\leq R$.
\end{proof}

\begin{theorem}\label{UE-gluing}
Let $X$ be a metric space. Suppose that for all $R>0$ and
$\varepsilon>0$ there is a partition of unity
$(\varphi_i)_{i\in I}$ on $X$ such that
\begin{ilist}
  \item  $\forall x,y\in X$ with $d(x,y)\leq R$,
     $\sum _{i\in I}|\varphi_i (x)-\varphi_i (y)|\leq
         \varepsilon.$
  \item $(\varphi_i)_{i\in I}$ is subordinated to an equi-uniformly
      embeddable cover $(U_i)_{i\in I}$
      of $X$.
\end{ilist}
Then $X$ is  uniformly embeddable.
\end{theorem}
\begin{proof}
Let $R>0$ and $\epsilon >0$ be given.
We construct a Hilbert space valued function $\eta$ on $X$ satisfying
the conditions in Proposition~\ref{char-UE}.
By assumption there is a
cover $\U =(U_i)_{i\in I}$ of $X$ which is equi-uniformly
embeddable and there is a partition of unity $(\varphi_i)_{i\in
I}$ subordinated to $\U$ such that $\forall x,y\in X$ with
$d(x,y)\leq R$,
\begin{equation*}
  \sum _{i\in I}|\varphi_i (x)-\varphi_i (y)|\leq
      \varepsilon^2/4.
\end{equation*}
Since $\U$ is equi-uniformly embeddable so is $\U_R=\{U_i(R):i \in
I\}$ where as above $U_i(R)=\{x\in X:d(x,U_i)\leq R\}$.  Therefore
 there exist Hilbert space valued maps
$\xi_{i}:U_i(R) \to \mathcal{H}_{i}$, with $\|\xi_{i}(x)\|=1$ for
all $x \in U_i(R)$ and such that
\begin{ilist}{\setcounter{ictr}{3}}
  \item $  \sup \{ \|\xi_{i}(x)-\xi_{i}(y)\|: d(x,y) \leq R,\, x,y\in U_i(R)\} \leq
       \varepsilon /2 $ for all $i\in I $.
  \item $\lim_{S \to \infty}
  \sup_{i\in I }\sup \{ |\langle\xi_{i}(x),\xi_{i}(y)\rangle|:
      d(x,y) \geq S,\, x,y\in U_i(R)\} =0 .$
\end{ilist}
We extend each $\xi_{i}$ to $X$ by setting $\xi_{i}(x)=0$ for $x
\in X\setminus U_i(R)$. Define $\eta:X \to \mathcal{H}=\oplus_{i
\in I} \mathcal{H}_{i}$, $\eta(x)=(\eta_i(x))_{i\in I}$ by
\begin{equation}\label{eta}
  \eta_i(x)= \varphi_i(x)^{1/2}\xi_{i}(x).
\end{equation}
One verifies that $\|\eta(x)\|=1$, $\forall\, x \in X$. Let $x,y
\in X$ with $d(x,y)\leq R$. Consider $\alpha(x,y),\beta(x,y)\in
\mathcal{H}$ with components
\[\alpha_i(x,y)=\varphi_i(x)^{1/2}(\xi_{i}(x)-\xi_{i}(y)),\]
\[\beta_i(x,y)=(\varphi_i(x)^{1/2}-\varphi_i(y)^{1/2})\xi_{i}(y). \]
Note that $\alpha(x,y)$ and $\beta(x,y)$ are well-defined because
of the following
 norm estimates.
\[\|\alpha(x,y)\|^2=\sum_{i\in
I}\varphi_i(x)\|\xi_{i}(x)-\xi_{i}(y)\|^2,
\]
where the summation is done for those $i$ with $x\in U_i$. If
$d(x,y)\leq R$ and $x \in U_i$ then $y \in U_i(R)$, so that using
(iv) we obtain $\|\alpha(x,y)\|\leq \varepsilon/2$. Since
$|a^{1/2}-b^{1/2}|^{ 2}\leq |a-b|$ we have
\begin{gather*}
 \|\beta(x,y)\|^2=\sum_{i\in I}\,\| \left
(\varphi_i(x)^{1/2}-\varphi_i(y)^{1/2})\,\xi_{i}(y)\right \|^2
 \leq
     \sum_{i\in I }|\varphi_i(x)^{1/2}-\varphi_i(y)^{1/2}|^2 \\
\leq \sum_{i\in I }|\varphi_i(x) -\varphi_i(y) |
\leq\varepsilon^2/4,
\end{gather*}
hence $\|\beta(x,y)\|\leq \varepsilon/2$. Therefore
\begin{equation*}
  \|\eta(x)-\eta(y) \|  =\|\alpha(x,y)+\beta(x,y)\|\leq
  \|\alpha(x,y)\|+\|\beta(x,y)\|\leq \varepsilon
\end{equation*}
 whenever $d(x,y)
\leq R$. In order to prove the support condition (ii) of
Proposition~\ref{char-UE} we use  that $\varphi_{i}$ vanishes
outside $U_i$ and the Cauchy-Schwarz inequality. Thus for any $x,y
\in X$ with $d(x,y)\geq S$ we have:
\begin{gather*}
      |\langle\;\eta(x),\eta(y)\rangle| \leq
       \sum_{i\in I}\varphi_i(x)^{1/2} \varphi_i(y)^{1/2}
      |\langle\;\xi_{i}(x),\xi_{i}(y)\rangle|\leq\\
       \sup_{i\in I }\sup \{ |\langle\xi_{i}(x'),\xi_{i}(y')\rangle|:
      d(x',y') \geq S,\, x',y'\in U_i\}.
\end{gather*}
In view of (v), this concludes the proof.
\end{proof}

\begin{corollary}\label{fibrewise-general}
Let $p:X \to Y$ be a map of metric spaces with the property that
$\forall R>0$ $\exists S>0$ such that $d(p(x),p(x')) \leq S$
whenever $d(x,x')\leq R$. Suppose that $Y$ is exact. If for each
uniformly bounded cover $(U_i)_{i\in I}$ of $Y$, the family
$(p^{-1}(U_i))_{i\in I}$ of subspaces of $X$ is equi-uniformly
embeddable \textup{(}respectively, equi-exact\textup{)}, then $X$ is
\UE\ \textup{(}respectively, exact\textup{)}.
\end{corollary}
\begin{proof}
Let $R>0$ and $\varepsilon >0$ be given and let $S>0$ be as in the
statement.  Since $Y$ is exact, we find a uniformly bounded cover
$(U_i)_{i\in I}$ of $Y$,   together with  a partition of unity
$(\varphi_i)_{i\in I}$ as in Definition~\ref{def-exactness} with $S$
playing the role of $R$. Then $(\varphi_i\circ p)_{i\in I}$ is a
partition of unity on $X$ subordinated to $(p^{-1}(U_i))_{i\in I}$
and satisfying the assumptions of Theorem~\ref{UE-gluing}
\textup{(}respectively, ~\ref{exact-gluing}\textup{)}.
\end{proof}


\begin{corollary}\label{equivariant}
  Let $p:X \to Y$ be a Lipschitz map of metric spaces.  Assume
  that a group $\Gamma$ acts by isometries on both $X$ and $Y$, that
  the action on $Y$ is transitive and that $p$ is
  $\Gamma$-equivariant.  Assume $Y$ is exact.  If there exists
  $y_0\in Y$ such that for every $n\in \mathbb{N}$ the inverse image
  $p^{-1}(B(y_0,n))$ is uniformly embeddable
  \textup{(}respectively, exact\textup{)}
  then $X$ is uniformly embeddable
  \textup{(}respectively, exact\textup{)}.
\end{corollary}



\begin{proof}
  This follows immediately from Corollary~\ref{fibrewise-general}.
  Indeed, let $(U_i)_{i\in I}$ be a uniformly bounded cover of $Y$.
  Since $\Gamma$ acts transitively and isometrically on $Y$, there
  are elements $s_i \in \Gamma$ and $n\in \mathbb{N}$ such that
  $s_iU_i \subseteq B(y_0,n)$ for all $i\in I$.
  Since $s_ip^{-1}(U_i)=p^{-1}(s_iU_i)\subset p^{-1}(B(y_0,n))$,
  we see that the family $(p^{-1}(U_i))_{i\in I}$ is isometric
  to a family of subspaces of $p^{-1}(B(y_0,n))$. Since the latter
  space is uniformly embeddable (exact), we conclude
  that the family
  $(p^{-1}(U_i))_{i\in I}$  is equi-uniformly embeddable
  (equi-exact).
\end{proof}

\section{Partitions of unity coming from the combinatorics of asymptotic
dimension}

Let $\U$ be a cover of $X$.  A {\it Lebesgue number\/} for $\U$ is
a number $L>0$ with the property that any subset $B\subset X$ of
diameter less than $L$ is contained in some $U \in \U$.  A cover
$\U$ of $X$ has {\it multiplicity at most $k$\/} if any $x \in X $
belongs to at most $k$ members of $\U$. One way to construct
partitions of unity with Lipschitz properties is given by the
following proposition.

\begin{proposition}\label{partunit}
Let $\U$ be a cover of a metric space $X$ with multiplicity at
most $k+1$, \textup{(}$k \geq 0$\textup{)} and Lebesgue number
$L>0$.  For $U \in \U$ define
\begin{equation*}
    \varphi_U(x) =  \frac{d(x,X\setminus U)}{\sum_{V \in \U} d(x,X\setminus V)}.
\end{equation*}
Then $(\varphi_U)_{U \in \U}$ is a partition of unity on $X$
subordinated to the cover $\U$.  Moreover each $\varphi_U$
satisfies
\begin{equation}\label{lip1}
  |\varphi_U(x)-\varphi_U(y)|\leq \frac{2k+3}{L}\,d(x,y),
           \quad \forall x,y \in X,
\end{equation}
and the family $(\varphi_U)_{U\in \U}$ satisfies
\begin{equation}\label{lip2}
  \sum_{U\in \U}|\varphi_U(x)-\varphi_U(y)|\leq
      \frac{(2k+2)(2k+3)}{L}\,d(x,y), \quad \forall x,y \in X.
\end{equation}
\end{proposition}
\begin{proof}
This is folklore. See Bell's paper for a proof of \eqref{lip1}.
Note that \eqref{lip2} follows from \eqref{lip1} since any point
in $X$ belongs to at most $k+1$ distinct elements of the cover
$\U$.
\end{proof}




A metric space $X$ has \emph{asymptotic dimension} $\leq k$ if for
every $R>0$ there exists a uniformly bounded cover $\U$ of $X$
such that every ball of radius $R$ in $X$ meets at most $k+1$
elements of $\U$.

In the context of non-uniformly bounded covers we require several
closely related properties.  Let $\mathcal{V}$ be a family of
nonempty subsets of the metric space $X$.  The {\it
multiplicity\/} of $\mathcal{V}$ is the maximum number of elements
of $\mathcal{V}$ that contain a common point of $X$; the {\it
$R$-multiplicity\/} of $\mathcal{V}$ is the maximum number of
elements of $\mathcal{V}$ that meet a common ball of radius $R$ in
$X$.  If $d(U,V)> L$ for all $U,V \in \mathcal{V}$ with $U \neq V$
then $\mathcal{V}$  is {\it $L$-separated\/} ($L>0$).  Note that
if $\mathcal{V}$ consists of just one set then $\V$ is
$L$-separated (vacuously) for every $L>0$. A cover $\mathcal{U}$
of $X$ is {\it $(k,L)$-separated\/} ($k\geq 0$ and $L>0$) if there
is a partition of $\U$ into $k+1$ families
\begin{equation*}
\U=\U_0\cup\cdots\cup\U_k
\end{equation*}
such that each family $\U_i$ is $L$-separated. In particular $\U$
has multiplicity at most $k+1$.

We make two observations, which we will apply to covers that are
not necessarily uniformly bounded.  First, a $(k,2R)$-separated
cover has $R$-multiplicity $\leq k+1$.  Second, if $\U$ is a cover
of $X$ with $L$-multiplicity $\leq k+1$ then $L$ is a Lebesgue
number for the cover $\U_L$ obtained by enlarging the sets in $\U$
by taking their $L$-neighborhoods;
\begin{equation*}
  \U_L = \{\, U(L )\colon U\in \U \,\}, \quad
       U(L) = \{\, x\in X \colon d(x,U) \leq L \,\}.
\end{equation*}
Further, in this case, the cover $\U_L$ has multiplicity $\leq
k+1$.

We summarize the previous discussion in the following form.

\begin{lemma}\label{discussion}
A $(k,2L)$-separated cover of a metric space has $L$-multiplicity
$\leq k+1$.  If a cover $\U$ of a metric space has
$L$-multiplicity $\leq k+1$ then the enlarged cover $\U_L$ has
multiplicity $\leq k+1$ and Lebesgue number $L$. \qed
\end{lemma}





The following result was proven by Higson and Roe in the case of
discrete bounded geometry metric spaces \cite[Lemma~4.3]{Hig-Roe};
in our more general setting it follows immediately from
Proposition~\ref{partunit} and Definition~\ref{def-exactness}.

\begin{proposition}\label{asdim->exact}
A metric space of finite asymptotic dimension is exact. \qed
\end{proposition}




We now prove a natural generalization of this result, where
uniform boundedness of the cover is replaced by the appropriate
uniform versions of uniform embeddability and Property $A$,
defined earlier. We also provide the proper setting to generalize
the `union theorems' of Bell and Dranishnikov
\cite{Bell,Bell-Dran}.

\begin{theorem}\label{gluing}
Let $X$ be a metric space. Assume that for every $\delta >0$ there
is a $(k,L)$-separated cover $\U$ of $X$ with $k^2+1\leq L\delta$
and such that the family $\U$ is \UUE \textup{(}where each $U \in
\U$ is given the induced metric from $X$\textup{)}. Then $X$ is
\UE.  If instead we assume that the family $\U$ is equi-exact then
$X$ is exact.
\end{theorem}

\begin{proof}
The statement concerning uniform embeddability follows from
Lemma~\ref{discussion}, Proposition~\ref{partunit} and
Theorem~\ref{UE-gluing}; for exactness use
Theorem~\ref{exact-gluing} instead of \ref{UE-gluing}.


More precisely, given $R>0$ and $\varepsilon>0$ we fix  a number
$\delta$, $0<\delta<1/20R$. Then
\begin{equation*}
   k^2+1\geq 2 (2k+2)(2k+3)R\delta,
\end{equation*}
for all integers $k \geq 0$. By assumption there is a
$(k,2L)$-separated cover $\U$ of $X$ such that $\U$ is \UUE\ and
$k^2+1\leq 2 L \delta\varepsilon$. By Lemma~\ref{discussion} the
cover $\U_{L}$ has multiplicity $\leq k+1$ and Lebesgue number $L$.
Proposition~\ref{partunit} provides a partition of unity
subordinated to $\U_{L}$ with the following property:  for all $x$,
$y\in X$, if $d(x,y)\leq R$ then
\begin{equation*}
     \sum_{U(L)\in \U_L} |\varphi_{U(L)}(x)-\varphi_{U(L)}(y)|
        \leq \frac{(2k+2)(2k+3) R}{L} \leq \frac{k^2+1}{2L\delta}
        \leq \varepsilon.
\end{equation*}
Since the cover $\U$ is \UUE\ so is the cover $\U_L$. We conclude
the proof by applying Theorem~\ref{UE-gluing}.
\end{proof}




\begin{corollary}\label{octopus}
If a metric space $X$ is a union of finitely many \UE\ subspaces
then $X$ is \UE. A similar result is true for exact spaces.
\end{corollary}

\begin{proof}
By assumption there is a finite cover $\U$ of $X$ with each $U \in
\U$ \UE. Let $\delta$ be given. Let $k+1$ denote the cardinality
of $\U$ and choose $L$ such that $k^2+1\leq L\delta$.  Then  $\U$
is a $(k,L)$-separated cover of $X$ and $\U$ is equi-uniformly
embeddable. The result follows now from Theorem~\ref{gluing}.
\end{proof}


\begin{corollary}\label{infiniteoctopus}
  If a metric space $X$ is the union of an equi-uniformly embeddable
  family of subspaces $\U$ with the property that for every $L>0$
  there is a uniformly embeddable subspace $Y\subset X$ such that the
  family $\{\, U\setminus Y \colon U\in\U \,\}$ is
  $L$-separated then $X$ is uniformly embeddable.
  A similar result is true for exact spaces.
\end{corollary}
\begin{proof} Given $\delta>0$ we fix  $L\geq 2/\delta$.
Let $\U$ and $Y$  (depending on $L$) be as in the statement. We
apply Theorem~\ref{gluing} using the $(1,L)$-separated cover of
$X$ given by the families of sets $\U_0\cup \U_1$, where
$\U_0=\{\, Y \,\}$ and $\U_1= \{\, U\setminus Y \colon U\in\U
\,\}$. 
\end{proof}


\begin{corollary}\label{equivariant-asdim}
Let $p:X \to Y$ be a Lipschitz map of metric spaces. Assume that a
group $\Gamma$ acts by isometries on both $X$ and $Y$, that the
action on $Y$ is transitive and that $p$ is $\Gamma$-equivariant.
Assume $Y$ has finite asymptotic dimension. If there exists
$y_0\in Y$ such that for every $n\in \mathbb{N}$ the inverse image
$p^{-1}(B(y_0,n))$ is uniformly embeddable, then $X$ is uniformly
embeddable.
\end{corollary}
\begin{proof}
This follows from Proposition ~\ref{asdim->exact} in conjunction
with Corollary~\ref{equivariant}
\end{proof}

\begin{remark}
In the previous corollary, if we assume instead that
$p^{-1}(B_Y(y_0,n))$ is exact then we conclude that $X$ is exact.
The result so obtained is closely related to a result of Bell
concerning Property $A$ for groups acting on spaces of finite
asymptotic dimension (compare \cite[Theorem~1]{Bell}).
\end{remark}


\section{Relatively hyperbolic groups}

Let $\Gamma$ be a finitely generated group which is hyperbolic
relative to a finite family $\{H_1,\dots,H_n\}$ of subgroups.  We
prove that $\Gamma$ is uniformly embeddable if
and only if each subgroup $H_i$ is uniformly embeddable.
There are two analogous results in the literature:
Osin proved an analogous statement for finite asymptotic dimension
\cite{Osin} and Ozawa proved an analogous statement for exactness
\cite{Ozawa}.  We rely heavily on Osin's method (Ozawa's method is
completely different), and are indebted to Ozawa for alerting us
to Osin's paper.


If $A$ is a symmetric set of generators of $\Gamma$, we denote by
$d_A$ the corresponding left-invariant metric on $\Gamma$. If $B$
is another such set with $A\subset B$ the identity map
$p:(\Gamma, d_A)\to (\Gamma, d_B)$ is equivariant and
$d_B(p(x),p(y))\leq d_A(x,y)$.  Let $S$ be a a finite symmetric
set generating $\Gamma$.  Let
\begin{equation*}
\mathcal{H}=\bigcup_{k}\,(H_k\setminus{e})
\end{equation*}
Let $d_{S}$ and $d_{S\cup\mathcal{H}}$ be the left invariant
metrics on $\Gamma$ induced by $S$ and $S\cup\mathcal{H}$,
respectively.  For $n \geq 1$, let
\begin{equation*}
B(n)=\{g \in \Gamma:d_{S\cup\mathcal{H}}(g,e)\leq n\}.
\end{equation*}
{\it We always view $B(n)$ as a subspace of $\Gamma$ equipped with
the
  metric $d_S$\/}.  The following useful recursive decomposition of
$B(n)$ is contained in the proof of \cite[Lemma 12]{Osin}:
\begin{equation}\label{B1}
   B(1)=S\cup \bigl( \bigcup_k H_k \bigr)
\end{equation}
\begin{equation}\label{Bn}
   B(n)=\bigl(\bigcup_k B(n-1)H_k\bigr) \cup
           \bigl(\bigcup_{x \in S} B(n-1)x\bigr)
\end{equation}
\begin{equation}\label{Bn-1}
    B(n-1)H_k=\bigsqcup_{g \in R(n-1)} gH_k,
\end{equation}
where the final equality represents a partition of $B(n-1)H_k$
into disjoint cosets according to a fixed set of coset
representatives, $R(n-1)\subset B(n-1)$.


\begin{proposition}[Osin]\label{Osin'}
For every $L>0$ there exists $\kappa(L)>0$ such that if
\begin{equation*}
Y=\{x\in \Gamma \,\colon\,d_S(x,B(n-1))\leq \kappa(L)\,\}
\end{equation*}
then for each $k$
\begin{equation}\label{16'}
    B(n-1)H_k \subset
        Y\cup \bigl(\bigcup_{g\in R(n-1)} gH_k\setminus Y \bigr)
\end{equation}
and the subspaces $gH_k\setminus Y$, $g\in R(n-1)$ of
$(\Gamma,d_S)$ are $L$-separated.
\end{proposition}
\begin{proof} The statement is implicit in the proof of \cite[Lemma 12]{Osin}.
\end{proof}

\begin{proposition}\label{coarse-stab=UE}
    If each $H_k$ is uniformly embeddable so is
    $B(n)$.  A similar statement is true for exactness.
\end{proposition}
\begin{proof}
The proof is by induction.  For the basis, observe that $B(1)$ is
uniformly embeddable by \eqref{B1} and the finite union theorem
(Corollary~\ref{octopus}).  For the induction step, assume that
$B(n-1)$ is uniformly embeddable. Using again the finite union
theorem and \eqref{Bn} we are reduced to verifying that each
$B(n-1)H_k$ is uniformly embeddable.  This follows from the
infinite union theorem (Corollary~\ref{infiniteoctopus}) and
Proposition~\ref{Osin'}.

The proof for exactness is analogous (compare \cite{Bell}).
\end{proof}

Osin also proves the following result \cite{Osin} (although not
explicitly stated, the result is the content of Lemmas~$17$, $18$
and $19$) :

\begin{proposition}[Osin]\label{Osin-asdim0}
    The metric space $(\Gamma, d_{S\cup\mathcal{H}})$ has finite
    asymptotic dimension.  \qed
\end{proposition}

\begin{theorem}\label{UE-thm}
Let $\Gamma$ be a finitely generated group which is hyperbolic
relative to a finite family $\{H_1,\dots,H_n\}$ of subgroups. Then
$\Gamma$ is uniformly embeddable in a Hilbert space if and only if
each subgroup $H_i$ is uniformly embeddable in a Hilbert space.
\end{theorem}
\begin{proof}
If $\Gamma$ is uniformly embeddable, then so are its subgroups,
the $H_k$.  For the converse we apply
Corollary~\ref{equivariant-asdim} to the isometric actions of
$\Gamma$ on the metric spaces $X=(\Gamma, d_{S})$, $Y=(\Gamma,
d_{S\cup\mathcal{H}})$, where $p$ is the identity map and $y_0=e$.
Then $B(n)=p^{-1}(B_Y(e,n))$ which is uniformly embeddable by
Proposition~\ref{coarse-stab=UE}.
\end{proof}
It is clear that the analogous result for  C*-exact groups, due to
Ozawa \cite{Ozawa}, can be recovered arguing as above.

\bibliographystyle{amsplain}


\providecommand{\bysame}{\leavevmode\hbox
to3em{\hrulefill}\thinspace}
\providecommand{\MR}{\relax\ifhmode\unskip\space\fi MR }
\providecommand{\MRhref}[2]{%
  \href{http://www.ams.org/mathscinet-getitem?mr=#1}{#2}
} \providecommand{\href}[2]{#2}

\end{document}